\documentclass[12pt]{amsart}

\usepackage{amssymb}
\usepackage{amsthm}

\newcommand{\pmodulo}[1]{\;(\mathrm{mod}\;#1)}

\usepackage{fullpage, txfonts, exscale}

\theoremstyle{plain}
\newtheorem*{thm}{Theorem}
\theoremstyle{definition}
\newtheorem*{rem}{Remark}

\begin{document}

\title[Sums of three squares]{On sums of three squares}
\author[S. Choi]{S.K.K. Choi}
\address{Department of Mathematics\\
Simon Fraser University \\
Burnaby, British Columbia \\
Canada V5A 1S6}
\email{kkchoi@cecm.sfu.ca}
\thanks{Research of Stephen Choi was supported by NSERC of Canada.}
\author[A. Kumchev]{A.V. Kumchev}
\address{Department of Mathematics \\
1 University Station C1200 \\
The University of Texas at Austin \\
Austin, TX 78712 \\
U.S.A.}
\email{kumchev@math.utexas.edu}
\author[R. Osburn]{R. Osburn}
\address{Max-Planck Institut f{\"{u}}r Mathematik, Vivatsgasse 7, 53111 Bonn,
Germany}
\email{osburn@mpim-bonn.mpg.de}
 
\date{\today}


\begin{abstract}
Let $r_{3}(n)$ be the number of representations of a positive integer $n$ as a sum of
three squares of integers. We give two alternative proofs of a conjecture of Wagon concerning
the asymptotic value of the mean square of $r_3(n)$.
\end{abstract}

\maketitle
\section{Introduction}

Problems concerning sums of three squares have a rich history. It is a classical result of Gauss that
\[
  n = x_{1}^{2} + x_{2}^{2} + x_{3}^{2}
\]
has a solution in integers if and only if $n$ is not of the form $4^a(8k+7)$ with $a$, $k \in \mathbb Z$. Let $r_{3}(n)$ be the number of representations of $n$ as a sum of three squares (counting signs and order). It was conjectured by Hardy and proved by Bateman \cite{Ba} that
\begin{equation}\label{i.1}
  r_3(n) = 4\pi n^{1/2} \mathfrak S_3(n),
\end{equation}
where the singular series $\mathfrak S_3(n)$ is given by \eqref{SS3} with $Q = \infty$. 

While in principle this exact formula can be used to answer almost any question concerning $r_3(n)$, the ensuing calculations can be tricky because of the slow convergence of the singular series $\mathfrak S_3(n)$. Thus, one often sidesteps \eqref{i.1} and attacks problems involving $r_3(n)$ directly. For example, concerning the mean value of $r_3(n)$, one can adapt the method of solution of the circle problem to obtain the following
\[
  \sum_{n \le x} r_{3}(n) \sim \frac{4}{3}\pi x^{3/2}.
\]
Moreover, such a direct approach enables us to bound the error term in this asymptotic formula. An application of a result of Landau \cite[pp. 200--218]{landau} yields  
\[ 
  \sum_{n \le x} {r_{3}(n)} = \frac {4}{3}{\pi}x^{3/2} + O\big( x^{3/4 + \epsilon} \big)
\]
for all $\epsilon > 0$, and subsequent improvements on the error term have been obtained by Vinogradov
\cite{vin}, Chamizo and Iwaniec \cite{ci}, and Heath-Brown \cite{hb}.

In this note we consider the mean square of $r_3(n)$. The following asymptotic formula was conjectured by 
Wagon and proved by Crandall (see \cite{CW} or \cite{BC}).

\begin{thm} 
  Let $r_{3}(n)$ be the number of representations of a positive integer $n$ as a sum of three squares of integers. Then 
  \begin{equation}\label{i.2}
    \sum_{n \le x} {r_{3}(n)}^2 \sim \frac {8{\pi}^{4}}{21\zeta(3)} x^2.
  \end{equation}
\end{thm}

Apparently, at the time they discussed this conjecture Crandall and Wagon were unaware of the earlier work of M\"uller \cite{mul1, mul2}. He obtained a more general result which, in a special case, gives
\[
  \sum_{n \le x} {r_{3}(n)}^2 = Bx^2 + O\big( x^{14/9} \big),
\]
where $B$ is a constant. However, since in M\"uller's work $B$ arises as a specialization of a more general (and more complicated) quantity, it is not immediately clear that $B = \frac 8{21}\pi^4/\zeta(3)$. The purpose of this paper is to give two distinct proofs of this fact: one that evaluates $B$ in the form given by M\"uller and a direct proof using the Hardy--Littlewood circle method.

\section{A direct proof: the circle method}

Our first proof exploits the observation that the left side of \eqref{i.2} counts solutions of the equation
\[
  m_1^2 + m_2^2 + m_3^2 = m_4^2 + m_5^2 + m_6^2
\]
in integers $m_1, \dots, m_6$ with $|m_j| \le x$. This is exactly the kind of problem that the circle method was designed for. The additional constraint $m_1^2 + m_2^2 + m_3^2 \le x$ causes some technical difficulties, but those are minor.

Set $N = \sqrt x$ and define
\[
  f(\alpha ) = \sum_{m \le N} e\big( \alpha m^2 \big),
\]
where $e(z) = e^{2\pi iz}$. Then for an integer $n \le x$, the number $r^*(n)$ of representations of $n$ as a sum of three squares of \emph{positive} integers is
\[
  r^*(n) =  \int_0^1 f(\alpha )^3 e(-\alpha n) \, d\alpha.
\]
Since $r_3(n) = 8r^*(n) + O(r_2(n))$, where $r_2(n)$ is the number of representations of $n$ as a sum of two squares, we have
\begin{equation}\label{2.0}
  \sum_{n \le x} r_3(n)^2 = 64\sum_{n \le x} r^*(n)^2 + O\big( x^{3/2 + \epsilon} \big).
\end{equation}
Therefore, it suffices to evaluate the mean square of $r^*(n)$. Let
\[
  P = N/4 \qquad \text{and} \qquad Q = N^{1/2}.
\]
We introduce the sets
\[
  \mathfrak M (q, a) = \big\{ \alpha \in \big[ Q^{-1}, 1 + Q^{-1} \big] : |q\alpha - a| \le PN^{-2} \big\}
\]
and 
\[
  \mathfrak M = \bigcup_{q \le Q}\bigcup_{\substack{1 \le a \le q \\ (a, q) = 1}} \mathfrak M (q, a), \qquad \mathfrak m = \big[ Q^{-1}, 1 + Q^{-1} \big] \setminus \mathfrak M.
\]
We have
\begin{equation}\label{2.1}
\begin{split}
  r^*(n) &= \bigg( \int_{\mathfrak M} + \int_{\mathfrak m} \bigg) f(\alpha )^3 e(-\alpha n) \, d\alpha \\
  &= r^*(n, \mathfrak M) + r^*(n, \mathfrak m), \qquad \text{say}. 
\end{split}
\end{equation}

We now proceed to approximate the mean square of $r^*(n)$ by that of $r^*(n, \mathfrak M)$. By \eqref{2.1} and Cauchy's inequality,
\begin{equation}\label{2.2}
  \sum_{n \le x} r^*(n)^2 = \sum_{n \le x} r^*(n, \mathfrak M)^2 + O\big( (\Sigma_1\Sigma_2)^{1/2} + \Sigma_2 \big),
\end{equation}
where 
\[
  \Sigma_1 = \sum_{n \le x} |r^*(n, \mathfrak M)|^2, \qquad 
  \Sigma_2 = \sum_{n \le x} |r^*(n, \mathfrak m)|^2.
\]

By Bessel's inequality,
\begin{equation}\label{2.3}
  |\Sigma_2| = \sum_{n \le x} \left| \int_{\mathfrak m} f(\alpha)^3 e(-\alpha n) \, d\alpha \right|^2 \le \int_{\mathfrak m} |f(\alpha)|^6 \, d\alpha. 
\end{equation}
By Dirichlet's theorem of diophantine approximation, we can write any real $\alpha$ as $\alpha = a/q + \beta $, where
\[
  1 \le q \le N^2P^{-1}, \qquad (a, q) = 1, \qquad |\beta| \le P/(qN^2). 
\]
When $\alpha \in {\mathfrak m}$, we have $q \ge Q$, and hence Weyl's inequality (see Vaughan \cite[Lemma 2.4]{Va}) yields
\begin{equation}\label{2.4}
  |f(\alpha )| \ll N^{1 + \epsilon}\big( q^{-1} + N^{-1} + qN^{-2} \big)^{1/2} \ll N^{1 + \epsilon} Q^{-1/2}. 
\end{equation}
Furthermore, we have
\begin{equation}\label{2.5}
  \int_0^1 |f(\alpha)|^4 \, d\alpha \ll N^{2 + \epsilon}, 
\end{equation}
because the integral on the left equals the number of solutions of
\[
  m_1^2 + m_2^2 = m_3^2 + m_4^2
\]
in integers $m_1, \dots, m_4 \le N$. For each choice of $m_1$ and $m_2$, this equation has $\ll N^{\epsilon}$ solutions. Combining \eqref{2.3}--\eqref{2.5} and replacing $\epsilon$ by $\epsilon/3$, we obtain
\begin{equation}\label{2.6}
  \Sigma_2 \ll N^{4 + \epsilon}Q^{-1}.
\end{equation}
Furthermore, another appeal to Bessel's inequality and appeals to \eqref{2.5} and to the trivial estimate $|f(\alpha)| \le N$ yield
\begin{equation}\label{2.7}
  \Sigma_1 \le \int_{\mathfrak M} |f(\alpha)|^6 \, d\alpha \le \int_0^1 |f(\alpha)|^6 \, d\alpha \ll N^{4 + \epsilon}.
\end{equation}

We now define a function $f^*$ on $\mathfrak M$ by setting
\[
  f^*(\alpha ) = q^{-1}S(q,a)v(\alpha - a/q) \qquad \text{for } \alpha \in \mathfrak M(q, a) \subseteq \mathfrak M; 
\]
here
\[
  S(q, a) = \sum_{1 \le h \le q} e \big( ah^2/q \big), \qquad
  v(\beta) = \frac 12 \sum_{m \le x} m^{-1/2} e(\beta m).
\]
Our next goal is to approximate the mean square of $r^*(n, \mathfrak M)$ by the mean square of the integral
\[
  R^*(n) = \int_{\mathfrak M} f^*(\alpha )^3 e(-\alpha n) \, d\alpha.
\]
Similarly to \eqref{2.2},
\begin{equation}\label{2.8}
  \sum_{n \le x} r^*(n, \mathfrak M)^2 = \sum_{n \le x} R^*(n)^2 + O \big( \Sigma_3 + (\Sigma_1\Sigma_3)^{1/2} \big),
\end{equation}
where 
\begin{equation}\label{2.9}
  \Sigma_3 = \sum_{n \le x} \bigg| \int_{\mathfrak M} \big[ f(\alpha)^3 - f^*(\alpha)^3 \big] e(-\alpha n) \, d\alpha \bigg|^2 \le \int_{\mathfrak M} \big| f(\alpha)^3 - f^*(\alpha)^3 \big|^2 \, d\alpha,
\end{equation}
after yet another appeal to Bessel's inequality. By \cite[Theorem 4.1]{Va}, when $\alpha \in \mathfrak M(q, a)$,
\[
  f(\alpha ) = f^*(\alpha) + O \big( q^{1/2 + \epsilon} \big).
\]
Thus,
\[
  \int_{\mathfrak M(q, a)} \big| f(\alpha)^3 - f^*(\alpha)^3 \big|^2 \, d\alpha
  \ll q^{1 + 2\epsilon} \int_{\mathfrak M(q, a)} \big( |f(\alpha)|^4 + q^{2 + 4\epsilon} \big) \, d\alpha,
\]
whence
\[
  \int_{\mathfrak M} \big| f(\alpha)^3 - f^*(\alpha)^3 \big|^2 \, d\alpha \ll
  Q^{1 + 2\epsilon} \int_0^1 |f(\alpha)|^4 \, d\alpha + PQ^{4 + 6\epsilon}N^{-2}.
\]
Bounding the last integral using \eqref{2.5} and substituting the ensuing estimate into \eqref{2.9}, we obtain
\begin{equation}\label{2.10}
  \Sigma_3 \ll QN^{2 + 2\epsilon} + PQ^4N^{-2 + 3\epsilon} \ll QN^{2 + 2\epsilon}.
\end{equation}
Combining \eqref{2.2}, \eqref{2.6}--\eqref{2.8}, and \eqref{2.10}, we deduce that
\begin{equation}\label{2.11}
  \sum_{n \le x} r^*(n)^2 = \sum_{n \le x} R^*(n)^2 + O \big( N^{4 + \epsilon}Q^{-1/2} + N^{3 + \epsilon} Q^{1/2} \big).
\end{equation}

We now proceed to evaluate the main term in \eqref{2.11}. We have
\[
  \int_{\mathfrak M(q, a)} f^*(\alpha)^3 e(-\alpha n) \, d\alpha
  = q^{-3}S(q, a)^3e(-an/q) \int_{\mathfrak M(q, 0)}v(\beta)^3e(-\beta n) \, d\beta,
\]
so
\[
  R^*(n) = \sum_{q \le Q} A(q, n)I(q, n),
\]
where
\[
  A(q, n) = \sum_{\substack{1 \le a \le q\\ (a, q) = 1}} q^{-3}S(q, a)^3 e(-an/q), \qquad
  I(q, n) = \int_{\mathfrak M(q, 0)}v(\beta)^3e(-\beta n) \, d\beta.
\]
Hence,
\begin{equation}\label{2.12}
  \sum_{n \le x} R^*(n)^2 = \sum_{n \le x} I(n)^2 \mathfrak S_3(n, Q)^2 + O\big( (\Sigma_4\Sigma_5)^{1/2} + \Sigma_5 \big),
\end{equation}
where
\begin{gather}\label{SS3}
  \mathfrak S_3(n, Q) = \sum_{q \le Q} A(q, n), \qquad I(n) = \int_{-1/2}^{1/2} v(\beta)^3 e(-\beta n) \, d\beta, \\
  \Sigma_4 = \sum_{n \le x} I(n)^2 \bigg( \sum_{q \le Q} |A(q, n)| \bigg)^2, \qquad
  \Sigma_5 = \sum_{n \le x} \bigg( \sum_{q \le Q} |A(q, n)(I(n) - I(q, n))| \bigg)^2. \notag
\end{gather}

By \cite[Theorem 2.3]{Va} and \cite[Theorem 4.2]{Va},
\begin{equation}\label{2.13}
  I(n) = \Gamma(3/2)^2\sqrt n + O(1) = \frac {\pi}4 \sqrt n + O(1), \qquad A(q, n) \ll q^{-1/2}.
\end{equation}
Furthermore, since $A(q, n)$ is multiplicative in $q$, \cite[Lemma 4.7]{Va} yields
\begin{align}\label{2.13a}
  \sum_{q \le Q} |A(q, n)| &\le \prod_{p \le Q} \big( 1 + |A(p, n)| + |A(p^2, n)| + \cdots \big) \\
  &\ll \prod_{p \le Q} \big( 1 + c_1(p, n)p^{-3/2} + 3c_1p^{-1} \big) \ll (nQ)^{\epsilon}, \notag
\end{align}
where $c_1 > 0$ is an absolute constant. In particular, we have
\begin{equation}\label{2.14}
  \Sigma_4 \ll N^{4 + \epsilon}.
\end{equation}

We now turn to the estimation of $\Sigma_5$. By Cauchy's inequality and the second bound in \eqref{2.13},
\[
  \Sigma_5 \ll (\log Q) \sum_{n \le x} \sum_{q \le Q} |I(n) - I(n, q)|^2 
\]
Another application of Bessel's inequality gives
\[
  \sum_{n \le x} |I(n) - I(n, q)|^2 \le 2\int_{P/qN^2}^{1/2} |v(\beta)|^6 \, d\beta.
\]
Using \cite[Lemma 2.8]{Va} to estimate the last integral, we deduce that
\[
  \Sigma_5 \ll \log Q \sum_{q \le Q} \big( q^2N^4P^{-2} + 1 \big) \ll N^2Q^{3 + \epsilon}.
\]
Substituting this inequality and \eqref{2.14} into \eqref{2.12}, we conclude that
\begin{equation}\label{2.15}
  \sum_{n \le x} R^*(n)^2 = \sum_{n \le x} I(n)^2\mathfrak S_3(n, Q)^2 + O \big( N^{3 + \epsilon}Q^{3/2} \big).
\end{equation}
We then use \eqref{2.13} and \eqref{2.13a} to replace $I(n)$ on the right side of \eqref{2.15} by $\frac {\pi}4\sqrt n$. We get 
\[
  \sum_{n \le x} I(n)^2\mathfrak S_3(n, Q)^2 = \frac {\pi^2}{16} \sum_{n \le x} n\mathfrak S_3(n, Q)^2 + 
  O \big( N^{3 + \epsilon} \big).
\]
Together with \eqref{2.11} and \eqref{2.15}, this leads to the asymptotic formula
\begin{equation}\label{2.16}
  \sum_{n \le x} r^*(n)^2 = \frac {\pi^2}{16} \sum_{n \le x} n\mathfrak S_3(n, Q)^2 + O \big( N^{4 + \epsilon}Q^{-1/2} + N^{3 + \epsilon}Q^{3/2} \big).
\end{equation}

Finally, we evaluate the sum on the right side of \eqref{2.16}. On observing that $\mathfrak S_3(n, Q)$ is in fact a real number, we have
\[
  \sum_{n \le t} \mathfrak S_3(n, Q)^2 = \sum_{q_1, q_2 \le Q} \sum_{ \substack{ 1 \le a_1 \le q_1\\ (a_1, q_1) = 1}} \sum_{ \substack{ 1 \le a_2 \le q_2\\ (a_2, q_2) = 1}} (q_1q_2)^{-3}S(q_1, a_1)^3 S(q_2, -a_2)^3 \sum_{n \le t} e\big( (a_1/q_1 - a_2/q_2)n \big).
\]
As the sum over $n$ equals $t + O(1)$ when $a_1 = a_2$ and $q_1 = q_2$ and $O(q_1q_2)$ otherwise, we get
\[
  \sum_{n \le t} \mathfrak S_3(n, Q)^2 = t\sum_{q \le Q} \sum_{ \substack{ 1 \le a \le q\\ (a, q) = 1}}
  q^{-6}|S(q, a)|^6 + O\big( \Sigma_6^2 \big),
\]
where 
\[
  \Sigma_6 = \sum_{q \le Q} \sum_{ \substack{ 1 \le a \le q\\ (a, q) = 1}} q^{-2}|S(q, a)|^3 \ll Q^{3/2}.
\]
We find that
\[
  \sum_{n \le t} \mathfrak S_3(n, Q)^2 = B_1t + O\big( tQ^{-1} + Q^3 \big),
\]
with
\[
  B_1 = \sum_{q = 1}^{\infty} \sum_{ \substack{ 1 \le a \le q\\ (a, q) = 1}} q^{-6}|S(q, a)|^6.
\]
Thus, by partial summation,
\[
  \sum_{n \le x} n\mathfrak S_3(n, Q)^2 = (B_1/2)x^2 + O\big( x^2Q^{-1} + xQ^3 \big).
\]
Combining this asymptotic formula with \eqref{2.16}, we deduce that
\[
  \sum_{n \le x} r^*(n)^2 = \frac {\pi^2}{32} B_1x^2 + O\big( x^{15/8 + \epsilon} \big).
\]
Recalling \eqref{2.0}, we see that \eqref{i.2} will follow if we show that
\[
  B_1 =  \frac{8\zeta (2)}{7\zeta (3)}.
\]
This, however, follows easily from the well-known formula (see \cite[\S 7.5]{Hua}) 
\begin{equation}\label{2.20}
  |S(q, a)| = \begin{cases}
    \sqrt q    & \text{if } q \equiv 1 \pmodulo 2, \\
    \sqrt {2q} & \text{if } q \equiv 0 \pmodulo 4, \\
    0          & \text{if } q \equiv 2 \pmodulo 4.
  \end{cases}
\end{equation}
Indeed, \eqref{2.20} yields
\[
  B_1 = \frac 43\sum_{q \text{ odd}} q^{-3}\phi(q) = \frac{8\zeta (2)}{7\zeta (3)},
\] 
where the last step uses the Euler product of $\zeta(s)$. This completes the proof of our theorem.

\section{Second Proof of Theorem}

Rankin \cite{Rankin} and Selberg \cite{Sel} independently
introduced an important method which allows one to study the
analytic behavior of the Dirichlet series

\begin{center}
$\displaystyle \sum_{n=1}^{\infty} \frac{a(n)}{n^s}$
\end{center}

\noindent where $a(n)$ are Fourier coefficients of a holomorphic
cusp form for some congruence subgroup of $\Gamma=SL_{2}(\mathbb
Z)$. Originally the method was for holomorphic cusp forms. Zagier
\cite{z} extended the method to cover forms that are not cuspidal
and may not decay rapidly at infinity. M{\"{u}}ller
\cite{mul1, mul2} considered the case where $a(n)$ is the
Fourier coefficient of non-holomorphic cusp or non-cusp form of
real weight with respect to a Fuchsian group of the first kind. It
is this last approach we wish to discuss. Note that if we apply a
Tauberian theorem to the above Dirichlet series, we then gain
information on the asymptotic behavior of the partial sum

\begin{center}
$\displaystyle \sum_{n\leq x} a(n)$.
\end{center}

We now discuss M{\"u}ller's elegant work. For details regarding discontinuous groups and
automorphic forms, see \cite{kubota, maass, mul1, rankin, roec1, roec2}. 
Let $\mathbb{H}=\{z\in \mathbb C: \Im(z)>0\}$ denote the upper half plane and $G=SL(2, \mathbb R)$ the special
linear group of all $2\times2$ matrices with determinant 1. $G$ acts on $\mathbb{H}$ by

\begin{center}
$\displaystyle z \mapsto gz=\frac{az+b}{cz+d}$
\end{center}

\noindent for $g=\left(\begin{matrix} a & b \\
c & d \\ \end{matrix} \right) \in G$. We write $y=y(z)=\Im(z)$. Thus we have

\begin{center}
$\displaystyle y(gz)=\frac{y}{|cz+d|^2}$.
\end{center}

\noindent Let $dx \hspace{.025in}dy$ denote the Lebesgue measure in the plane. Then the measure

\begin{center}
$\displaystyle d\mu=\frac{dx \hspace{.025in}dy}{y^2}$
\end{center}

\noindent is invariant under the action of $G$ on $\mathbb{H}$. A discrete subgroup $\Gamma$ of $G$ is called a
{\it Fuchsian group of the first kind} if its fundamental domain $\Gamma\backslash\mathbb{H}$ has finite volume. 
Let $\Gamma$ be a Fuchsian group of the first kind containing $\pm I$ where
$I$ is the identity matrix. Let $\mathcal{F}(\Gamma, \chi, k, \lambda)$ denote the space
of (non-holomorphic) automorphic forms of real weight $k$, eigenvalue
$\lambda=\frac{1}{4}-\rho^2$, $\Re(\rho) \geq 0$, and multiplier system $\chi$.
For $k \in \mathbb R$, $g \in SL(2,\mathbb R)$ and $f:\mathbb{H}\to \mathbb{C}$, we define the
stroke operator $|_{k}$ by

\begin{center}
$\displaystyle (f|_{k}g)(z):=\Bigg(\frac{cz+d}{|cz+d|}\Bigg)^{-k} f(gz)$
\end{center}

\noindent where $g=\left(\begin{matrix} a & b \\
c & d \\ \end{matrix} \right) \in \Gamma$. The transformation law for $f \in
\mathcal{F}(\Gamma, \chi, k, \lambda)$ is then

\begin{center}
$(f|_{k}g)(z)=\chi(g)f(z)$
\end{center}

\noindent for all $g \in \Gamma$. Automorphic forms $f \in \mathcal{F}(\Gamma, \chi, k, \lambda)$ have a Fourier expansion
at every cusp $\kappa$ of $\Gamma$, namely
\[
  A_{\kappa, 0}(y) + \sum_{n \neq 0} a_{\kappa, n}W_{(sgn \hspace{.025in} n)\frac{k}{2},
\rho}(4\pi|n+\mu_{\kappa}|y)e((n+\mu_{k})x),
\]
where $\mu_{\kappa}$ is the cusp parameter and $a_{\kappa, n}$ are the Fourier
coefficients of $f$ at $\kappa$. The functions $W_{\alpha,\rho}$ are Whittaker functions
(see \cite[\S 3]{mul1}), $A_{\kappa, 0}(y)=0$ if $\mu_{k} \neq 0$ and 

$$A_{\kappa,0}(y)=
\left \{ \begin{array}{l}
a_{\kappa,0}y^{1/2 + \rho} + b_{\kappa,0}y^{1/2 - \rho} \quad \mbox{if $\mu_{\kappa}=0$,
$\rho \neq 0$},\\
a_{\kappa,0}y^{1/2} + b_{\kappa,0}y^{1/2}\log y \quad  \mbox{if  $\mu_{\kappa}=0$,
$\rho =0$}.
\end{array}
\right. \\
$$

\noindent An automorphic form $f$ is called a cusp form if $a_{\kappa,0}=b_{\kappa,0}=0$
for all cusps $\kappa$ of $\Gamma$. Now consider the Dirichlet series
\[
  S_{\kappa}(f,s)=\sum_{n>0} \frac{|a_{\kappa, n}|^2}{(n+\mu_{\kappa})^s}.
\]
This series is absolutely convergent for $\Re(s)>2\Re(\rho)$ and has been shown
\cite{mul2} to have meromorphic continuation in the entire complex plane. In what
follows, we will only be interested in the case $f$ is not a cusp form.  
If $f$ is not a cusp form and $\Re(\rho)>0$, then $S_{\kappa}(f,s)$ has a simple pole at
$s=2\Re(\rho)$ with residue
\begin{equation} 
\beta_{\kappa}(f)=\underset{s=2\Re(\rho)}{\text{res}} S_{\kappa}(f,s) = (4\pi)^{2\Re(\rho)}
b^{+}(k/2, \rho)\sum_{\iota \in K} \varphi_{\kappa, \iota}(1+2\Re(\rho))|a_{\iota, 0}|^2,
\label{23}
\end{equation} 
where $K$ denotes a complete set of $\Gamma$-inequivalent cusps, $\varphi_{\kappa,
\iota}(1+2\Re(\rho))>0$ and $b^{+}(\frac{k}{2}, \rho)>0$  if $\rho+\frac{1}{2} \pm \frac{k}{2}$ is a
non-negative integer. For the definition of the functions $\varphi_{\kappa, \iota}$ and $b^{+}$, see Lemma
3.6 and (69) in \cite{mul2}. This result \eqref{23} and a Tauberian argument then provide the asymptotic behaviour of the summatory function 
\[
  \displaystyle \sum_{n \leq x} |a_{\kappa, n}|^2 |n+\mu_{\kappa}|^r. 
\]
 Precisely, we have (see \cite[Theorem 2.1]{mul1} or \cite[Theorem 5.2]{mul2}) that
\begin{equation}
\label{24}
\displaystyle \sum_{n \leq x} |a_{\kappa, n}|^2 |n+\mu_{\kappa}|^r = \sum_{z\in R}  \underset{s=z}{\text{res}} 
S_{\kappa}(f,s)\frac{x^{r+s}}{r+s} + O(x^{r+2\Re\rho -\gamma}(\log x)^g),
\end{equation}
 where $2\Re(\rho)+r \geq 0$, $R=\{\pm 2\Re(\rho), \pm 2i\Im(\rho), 0,
-r\}$, $\gamma=(2+8\Re(\rho))(5+16\Re(\rho))^{-1}$, and $g=$max$(0,b-1)$; $b$
denotes the order of the pole of $S_{\kappa}(f,s)(r+s)^{-1}x^{r+s}$ at
$s=2\Re(\rho)$ $(0 \leq b \leq 5)$. 

We now consider an application of \eqref{24}.
Let $Q \in \mathbb Z^{m \times m}$ be a non-singular symmetric
matrix with even diagonal entries and $q({\bf
x})=\frac{1}{2}Q[{\bf x}]=\frac{1}{2}{{\bf x}^T}Q{\bf x}$, ${\bf
x} \in \mathbb Z^{m}$, the associated quadratic form in $m\geq 3$
variables. Here we assume that $q({\bf x})$ is positive definite. 
Let $r(Q,n)$ denote the number of representations of
$n$ by the quadratic form $Q$. Now consider the theta function

\begin{center}
$\theta_{Q}(z)=\displaystyle \sum_{{\bf x} \in \mathbb Z^{m}}
e^{{\pi}izQ[{\bf x}]}$.
\end{center}

\noindent By \cite[Lemma 6.1]{mul1}, the Dirichlet series associated with the automorphic
form $\theta_{Q}$ is

\begin{center}
$(4\pi)^{-m/4}\zeta_{Q}(\frac{m}{4}+ s)$
\end{center}

\noindent where

\begin{center}
$\zeta_{Q}(s)=\displaystyle \sum_{n=1}^{\infty} \frac{r(Q,n)}{n^s}
= \sum_{{\bf x} \in \mathbb Z^{m}\setminus\{\bf{0}\}} q({\bf
x})^{-s}$
\end{center}

\noindent for $\Re(s)>m/2$. Using \eqref{24}, M{\"u}ller proved the following (see \cite[Theorem 6.1]{mul1})

\begin{thm}[M\"{u}ller] Let $q({\bf
x})=\frac{1}{2}Q[{\bf x}]=\frac{1}{2}{{\bf x}^T}Q{\bf x}$, ${\bf
x} \in \mathbb Z^{m}$ be a primitive positive definite quadratic form in $m\geq 3$ variables with integral coefficients. Then
\[ \sum_{n\leq x} r(Q,n)^2 = Bx^{m-1} + O\Big(x^{(m-1)\frac{4m-5}{4m-3}}\Big) \]
where 
\[ B=(4\pi)^{m/2}\frac{\beta_{\infty}(\theta_{Q})}{m-1} \]
and $\beta_{\infty}(\theta_{Q})$ is given by \eqref{23}.
\end{thm}

\noindent We are now in a position to prove our theorem in page 2.

\begin{proof}
 
We are interested in the case $q({\bf x})=x_{1}^{2}
+ x_{2}^{2} + x_{3}^{2}$ and so $r(Q,n)=r_3(n)$ counts the number of
representations of $n$ as a sum of three squares. By M\"{u}ller's Theorem above,

\begin{center}
$\displaystyle \sum_{n\leq x} r_3(n)^2 = Bx^2 + O\Big(x^{14/9}\Big)$
\end{center}

\noindent where $B$ is a computable constant. Specifically, we have by \eqref{23} (with $k=3/2$ and $\rho=1/4$)

\begin{center}
$\displaystyle
B=\frac{4\pi^2}{3-1}b^{+}(3/4, 1/4)\sum_{\iota \in K}
\varphi_{\infty, \iota}(3/2) |a_{\iota, 0}|^2$,
\end{center}

\noindent where $K$ denotes a complete set of
$\Gamma_0(4)$-inequivalent cusps  
and $a_{\iota, 0}$ is the $0$-th Fourier coefficient of $\theta_{Q}(z)$
at a rational cusp $\iota$. Choose $K=\{1, \frac{1}{2}, \frac{1}{4} \}$. Then by p. 145 and
(67) in \cite{mul1}, we have

\begin{center}
$\displaystyle |a_{\iota, 0}|^2=W_{\iota}^3 |G(S_{\iota})|^2$
\end{center}

\noindent where $\iota=u/w$, $(u,w)=1$, $w \geq 1$, $W_{\iota}$ is width of the cusp
$\iota$, and

\begin{center}
$\displaystyle |G(S_{\iota})|^2=2^{-3}w^{-3}\Big| \sum_{x=1}^{w} e(\frac{u}{w}x^2)
\Big|^{6}$.
\end{center}

As $W_{1/4}=W_{1/2}=1$, $W_{1}=4$, we have $|a_{1, 0}|^2=1$, $|a_{1/2, 0}|^2=0$,
and $|a_{1/4, 0}|^2=1$. An explicit description of the functions 
$\varphi_{\infty, \iota}(s)$ in the case $\Gamma_{0}(4)$ is given by (see (1.17)
and p. 247 in \cite{di})

\begin{center}
$\displaystyle \varphi_{\infty,
1/4}(s)=2^{1-4s}(1-2^{-2s})^{-1}{\pi}^{1/2}\frac{\Gamma(s-1/2)\zeta(2s-1)}{\Gamma(s)\zeta(2s)}$,
\end{center}

\begin{center}
$\displaystyle \varphi_{\infty,
1/2}(s)=\varphi_{\infty,1}(s)=2^{-2s}(1-2^{-2s})^{-1}(1-2^{1-2s}){\pi}^{1/2}\frac{\Gamma(s-1/2)\zeta(2s-1)}{\Gamma(s)\zeta(2s)}$.
\end{center} 

\noindent Thus for $s=3/2$, we have 

\begin{center}
$\displaystyle \varphi_{\infty,
1/4}(3/2)=2^{-5}(1-2^{-3})^{-1}{\pi^2}\frac{\zeta(2)}{\Gamma(3/2)\zeta(3)}$,
\end{center}

\begin{center}
$\displaystyle \varphi_{\infty,
1/2}(3/2)=\varphi_{\infty,1}(3/2)=2^{-3}(1-2^{-3})^{-1}(1-2^{-2}){\pi^2}\frac{\zeta(2)}{\Gamma(3/2)\zeta(3)}$.
\end{center}

\noindent Now, from p. 65 in \cite{mul2}, we have

\begin{center}
$\displaystyle b^{+}(3/4, 1/4)=G_{1/4, 1/4}^{*}(3/2)$.
\end{center}

\noindent By Lemma 3.3 and (16) in \cite{mul2},

\begin{center}
$\displaystyle G_{1/4, 1/4}^{*}(s)=\Gamma(s + 1/2)^{-1}$
\end{center}

\noindent and so $b^{+}(3/4, 1/4)={\Gamma(2)}^{-1}$. In total, 

$$
\begin{aligned}
B &=\frac{(4\pi)^2}{(3-1)}\frac{1}{\Gamma(2)} \Bigg(
2^{-3}(1-2^{-3})^{-1}(1-2^{-2}){\pi}^{1/2} \frac{\zeta(2)}{\Gamma(3/2)\zeta(3)} \\
&+ 2^{-5}(1-2^{-3})^{-1}{\pi}^{1/2}\frac{\zeta(2)}{\Gamma(3/2)\zeta(3)} \Bigg) \\
&= \frac{8{\pi}^{4}}{21\zeta(3)}. \\
\end{aligned}
$$

\noindent Thus 

\begin{center}
$\displaystyle\sum_{n \le x} {r_{3}(n)}^2 \sim
\frac{8{\pi}^{4}}{21\zeta(3)} x^2$.
\end{center}
\end{proof} 

\begin{rem} M\"{u}ller's Theorem can also be used to obtain the mean square
value of sums of $N > 3$ squares. Precisely, if $r_{N}(n)$ is the number of
representations of $n$ by $N > 3$ squares, then a calculation similar to the second proof
of our theorem yields (compare with Theorem 3.3 in \cite{BC})

\begin{center}
$\displaystyle \sum_{n\leq x} r_N(n)^2 = W_{N}x^{N-1} +
O\Big(x^{(N-1)\frac{4N-5}{4N-3}}\Big) $
\end{center}
where
\begin{center}
$\displaystyle
W_{N}=\frac{1}{(N-1)(1-2^{-N})}\frac{\pi^N}{\Gamma(N/2)^2}\frac{\zeta(N-1)}{\zeta(N)}$.
\end{center}
\end{rem}

\section*{Acknowledgments}
The authors would like to thank Wolfgang M{\"{u}}ller for his comments
regarding the second proof of the theorem. The second author would like
to take this opportunity to express his gratitude to the Mathematics Department
at the University of Texas at Austin for the support during the past three years. 
The third author would like to thank the Max-Planck-Institut f{\"u}r Mathematik 
for their hospitality and support during the preparation of this paper.

\end{document}